\documentclass[12pt]{amsart}
\usepackage{amssymb}

\newtheorem{theorem}{Theorem}[section]

\theoremstyle{definition}

\newtheorem{corollary}[theorem]{Corollary}
\newtheorem{proposition}[theorem]{Proposition}
\theoremstyle{remark}
\newtheorem{remark}[theorem]{Remark}

\begin{document}
\def\C{\mathbb C}
\def\R{\mathbb R}
\def\X{\mathbb X}
\def\cA{\mathcal A}
\def\cT{\mathcal T}
\def\Z{\mathbb Z}
\def\Y{\mathbb Y}
\def\Z{\mathbb Z}
\def\N{\mathbb N}
\def\cal{\mathcal}

\title[on the Solvability of Some Abstract Differential Equations]{Some Remarks on the Solvability of 
Some Abstract Differential Equations}
\author{Toka Diagana} 
\address{Department of Mathematics, Howard University, 2441 6th Street N.W. Wasington D.C. 20059, USA} 
\email{tdiagana@howard.edu}
\begin{abstract}
This paper is concerned with the solvability of some abstract differential
equation of type
$\dot u(t) + Au(t) + Bu(t)  \ni f(t),\; t \in (0,T], \; u(0)
= 0$,
where $A$ is a linear selfadjoint operator and $B$ is a nonlinear (possibly multi-valued)
maximal monotone operator in a real Hilbert space ${\mathbb H}$ with
the normalization $0 \in B (0)$.
We use the concept of variational sum introduced
by H. Attouch, J.-B. Baillon, and M. Th\'era, to investigate solutions to the given differential equation. 
\end{abstract}
\date{}
\maketitle

\section{Introduction} Our aim in this paper is to investigate on the solvability
of some abstract differential equation of type 

\[ (1.1) \; \; \;
\left \{ \begin{array}{cl}
\dot u(t) + Au(t) + Bu(t)  \ni f(t), \; \; t \in (0, T]\\
u(0) = 0
\end{array} \right.
\]
\noindent where $A$ is a linear selfadjoint monotone operator and $B$ (the
operator $B$ some times equals to $\partial
\psi$, the subdifferential of a convex lower semicontinous proper function
$\psi: \; {\mathbb H} \mapsto {\mathbb R} \cup {+ \infty}$) is a nonlinear
maximal
monotone 
operator (possibly multi-valued) in the real Hilbert space
$( {\mathbb H} \; ; \; \langle \; , \; \rangle )$ with the normalization $0 \in B (0)$. The function $f$ belongs to
$L^2( 0 , T ; {\mathbb H})$, where $L^2( 0 , T ; {\mathbb H})$ is the Hilbert space endowed
with the inner product 
$$ (1.2) \; \; \; \; \; \; \langle \langle u \; , \; v  \rangle \rangle = \int_{0}^{T} \; \langle u(t) \;
, \; v(t) \rangle
\; dt$$
It is convenient to write $(1.1)$ of the form
$$(1.3) \; \; \; \; \; \; {\mathcal S} u + {\mathcal A} u + {\mathcal B}
u \ni f$$
where ${\mathcal S}$ is defined in $L^2( 0 , T ; {\mathbb H})$ by
\[ \left \{ \begin{array}{cl}
D({\mathcal S}) = \{ u \in {\mathbb H^1} (0 , T ; {\mathbb H}) \; | \; u(0) = 0\}\\
{\mathcal S} u = \dot u, \; \; \; \forall u \in D({\mathcal S})
\end{array} \right.
\]
The operator ${\mathcal A}$ is defined in $L^2(0 , T ; {\mathbb H})$ by:
$ u \in D({\mathcal A}) \; \; \mbox{and} \; \; {\mathcal A}u = v$
iff $u, v  \in L^2(0 , T ; {\mathbb H})$ and $A u (.) = v(.)$. In the
same way, the operator ${\mathcal B}$ is defined in $L^2(0 , T; {\mathbb
H})$ by:
$ u \in D({\mathcal B}) \; \; \mbox{and} \; \; {\mathcal B}u = v$
iff $u, v  \in L^2(0 , T ; {\mathbb H})$ and ${\mathcal B} u (.) = v(.)$.

It is well-known that the corresponding operators ${\mathcal S}$, ${\mathcal A}$, and ${\mathcal B}$ are maximal monotone on $L^2(0 , T ; {\mathbb H})$, see, e.g., \cite{har}.

Recall that various types of equation (1.1) have been investigated in
the past decades by
several mathematicians, see, e.g., \cite{am, ben, br-gr-pa, br-st, cl-eg,
fu-ka, gi-mi, har, ka, ka1, fuj-kat}. 

This paper is concerned with the solvability to the equation (1.1).
Assuming that $A$ is a linear, selfadjoint monotone operator (possibly unbounded)
and that the operator $B$ (possibly multi-valued) is 
maximal accretive. Under suitable assumptions, we shall establish the
solvability of the equation (1.1). Recall that the innovating idea in
this paper is
the use of the powerful concept of the variational sum, introduced
by H. Attouch, J B. Baillon, and M. Th\'era in \cite{at-ba-th}.

In what follows,
we assume that $A$
is a linear, self-adjoint monotone operator (possibly unbounded) and that
$B$ is maximal monotone(possibly multi-valued), with the normalization
$0 \in B(0)$.  Instead of considering (1.1), we consider the solvability
of (1.3).

Recall that the Moreau-Yosida approximation of ${\mathcal B}$ is defined
as follows
$${\mathcal B}_{\lambda} = {\frac{1}{\lambda}} \; [ I - (I + \lambda {\mathcal
B})^{-1}] \; \;
\; \lambda \; > \; 0$$
The Moreau-Yosida approximation ${\mathcal B}_{\lambda}$, $\; \; \lambda \; > \; 0$ is
an everywhere defined, Lipschitz continuous, and maximal monotone operator,
see, e.g., \cite{br}. Also, note that ${\mathcal B}_{\lambda}$ is single-valued
and that ${\mathcal B}_{\lambda} \subset {\mathcal B} \; ( I + \lambda {\mathcal
B})^{-1}$ in the
sense of the corresponding graphs.

Now recall the definition and some details about the concept of the variational sum $({\mathcal A} + {\mathcal B})_{v}$ of 
${\mathcal A}$ and ${\mathcal B}$.
Let ${\mathcal{F}}$ be the {\it filter} of all the pointed
neighborhoods of the origin $(0,0)$ in the set
$I = \{ (\lambda, \mu) \in {\mathbb{R}}^2 : \; \;
\lambda, \; \mu \; \geq 0, \; \; \lambda + \mu \not= 0 \; \}$
and $\lim \inf_{{\mathcal{F}}}$ for $\lim_{\lambda \to 0, \; \mu \to 0, \;
(\lambda , \mu ) \in I}$.
The variational sum of the maximal monotone operators ${\mathcal A}$ and ${\mathcal B}$ is
defined as
$$ (1.4) \; \; \; \; \; \; ({\mathcal A} \; + \; {\mathcal B})_{v} = \lim \inf_{{\mathcal{F}}} \; ({\mathcal A}_{\lambda} + {\mathcal B}_{\lambda})$$
where the limit Inferior is understood in the sense of Kuratowski-Painlev\'e,
when ${\mathcal A}$ and ${\mathcal B}$ are identified with their graphs. The equation $(1.4)$ can
be equivalently expressed in terms of resolvents, that is,
for any $w \in L^2(0 , T ; {\mathbb H})$, the family $\{ u_{\lambda, \mu} : \;
\; (\lambda , \mu ) \in I \}$ of solutions of
$$(1.5) \; \; \; \; \; \; \; u_{\lambda , \mu} + {\mathcal A}_{\lambda} u_{\lambda , \mu} +
{\mathcal B}_{\mu} u_{\lambda , \mu} u_{\lambda , \mu} \; \ni w$$
converges (with respect to the {\it filter} ${\mathcal{F}}$) to
the solution $u$ of
$$(1.6) \; \; \; \; \; \; \; u + ({\mathcal A} \; + \; {\mathcal
B})_{v} u \; \ni w$$
More details about the concepts of the algebraic sum,
variational sum, generalized sum or the extended sum 
can be found in \cite{at, at-ba-th, br, br-gr-pa, cl-eg, da-gri, dia1, dia2,
nu, re-th, rev-the}. 
Nevertheless, recall that among the main motivations of the concept 
of variational sum, we have the fact that
the algebraic sum of two operators
in not always well-adapted to problems arising in mathematics, see \cite{dia2,
dia3, dia4}, for
examples.

\section{Existence of solutions} Consider the equation (1.3) 
in the Hilbert space $L^2(0 , T ; {\mathbb H})$. 
Thus, the existence problem of solutions to (1.3) is equivalent to finding conditions for which
the algebraic sum ${\mathcal A+B}$ is a maximal accretive operator in $L^2(0
, T ; {\mathbb H})$.
As stated in the introduction, the algebraic sum is not well-adapted to the
present situation. To overcome such a difficulty, we shall deal this the variational sum $({\mathcal A} + {\mathcal B})_{v}$ of ${\mathcal A}$ and ${\mathcal
B}$ and compute
it.

\subsection{Revolvent Commuting case}
Let ${\mathcal A}$ and ${\mathcal B}$ be the operators described above, where ${\mathcal
B}$ is supposed to be a linear operator. Assume that they commute
in the sense of resolvent, that is,
$$(2.1) \; \; \; \; \; ( I + \lambda {\mathcal A})^{-1} (I + \mu {\mathcal B})^{-1}
= (I + \mu {\mathcal B})^{-1} ( I + \lambda {\mathcal A})^{-1}, \; \;
\forall \lambda, \mu > 0$$ In this case, it is well-known that
the algebraic sum ${\mathcal A} + {\mathcal B}$ is closable, and by
a result of Da Prato and P. Grisvard (see \cite{da-gri}), we also know 
$\overline{{\mathcal A} + {\mathcal B}}$ is m-accretive.

We have
\begin{proposition}
Let ${\mathcal A}$, ${\mathcal B}$ be the corresponding operators
to $A$ and $B$ described above. Assume that ${\mathcal B}$ is linear and that
(2.1) holds; then the
problem (1.3) has a unique solution.
\end{proposition}

\begin{proof}
Since ${\mathcal A}$ (densely defined) and ${\mathcal B}$ are
m-accretive satisfying (2.1). Then according to \cite{da-gri}, 
$- \overline{{\mathcal A} + {\mathcal B}}$ is m-dissipative. Therefore
the Hille-Yosida theorem guarantees the existence
of a unique solution to (1.3). Since (1.1) and (1.3) are equivalent,
then so does (1.1).
\end{proof}
\begin{remark}
Considering the given linear m-accretive operators ${\mathcal A}$ and ${\mathcal B}$ described above. If one of them generates an analytic semigroup and that
(2.1) holds, it is also known that $- \overline{{\mathcal A} + {\mathcal B}}$ is m-dissipative. This is in fact a consequence of a result due to
Dore and Venni (see \cite{do-ve}).
\end{remark}

\subsection{Case where ${\mathcal B} = \partial \psi$}
We assume that the operator ${\mathcal B} = \partial \phi$ is the subdifferential of
a convex semicontinuous proper function $\psi: \; L^2(0 , T; {\mathbb H}) \mapsto {\mathbb
R} \cup {+ \infty}$. Under previous assumptions, it is well-known that
the Moreau-Yosida approximation ${\mathcal B_{\lambda}} = (\partial \psi)_{\lambda}
= \partial \psi_{\lambda}$ where
$$ \psi_{\lambda} (x) = \inf_{v \in {\mathbb H}} \{ \psi(v) + {\frac{1}{2
\lambda}} \| x - v \|^2 \}$$

Recall that since ${\mathcal A}$ is a self-adjoint monotone operator, then
it can be expressed as ${\mathcal A} = \partial \phi$, the subdifferential of the
convex semicontinuous proper functional $\phi: \; L^2(0 , T ; {\mathbb H}) \mapsto {\mathbb
R} \cup {+ \infty}$ defined by

\[ (2.2) \; \; \; \phi(u) = \; \;
\left \{ \begin{array}{cl}
\displaystyle{\frac{1}{2}} \| {\mathcal A^{\frac{1}{2}}} u \|^2 \; \; \; \mbox{if} 
\; u \in D({\mathcal A^{\frac{1}{2}}})\\
+ \infty \; \; \; \mbox{elsewhere}
\end{array} \right.
\]

Since ${\mathcal B}_{\lambda}$ is m-accretive and Lipschitz $\forall \lambda > 0$, it is well-known that ${\mathcal A} + {\mathcal B}_{\lambda}$
is maximal monotone operator, see, e.g.,\cite{br}. Therefore, for all $w \in L^2(0
, T ; {\mathbb H})$,
the equation 
$$(2.3) \; \; \; \; \; \; u_{\lambda} + {\mathcal A} u_{\lambda} + {\mathcal B}_{\lambda} u_{\lambda} = w, \; \; \lambda > 0$$
has a unique solution $u_{\lambda} \in D({\mathcal A})$. We also know that
there exists a unique $u \in D(({\mathcal A} + {\mathcal B})_{v})$ such that
$u_{\lambda}$ converges to $u$ as $\lambda$ goes to zero, and
$$(2.4) \; \; \; \; \; \; w \; \in \; u + ({\mathcal A} + {\mathcal B})_{v}
u $$

\begin{theorem}
{\it Let ${\mathcal A}$ and ${\mathcal B}$ be the operators described above
such that 
${\mathcal A} = \partial \phi$ and ${\mathcal B} = \partial \psi$
on $L^2( 0 , T; {\mathbb H})$. Assume that $D( \phi )
\cap D( \psi ) \not = \varnothing$. Then the variational sum of ${\mathcal
A}$ and ${\mathcal B}$ is given by
$$ ({\mathcal A} + {\mathcal B})_{v} = \partial (\phi + \psi)$$
}
\end{theorem}

\begin{proof}
This is a consequence of ( Theorem 7.2, \cite{at-ba-th} ), it is straightforward.
\end{proof}

\begin{corollary}
{\it Let ${\mathcal A}$ and ${\mathcal B}$ be the operators described above
such that 
${\mathcal A} = \partial \phi$ and ${\mathcal B} = \partial \psi$
on $L^2( 0 , T; {\mathbb H})$. Assume that $D( \phi )
\cap D( \psi ) \not = \varnothing$. Then the problem given by

$$(2.6) \; \; \; \; \; \; \dot u  + \partial( \phi + \psi) u  \ni f $$
has a unique solution.}
\end{corollary}

\begin{remark}
In the case where $\overline{{\mathcal A} + {\mathcal B}}$ is an m-accretive
operator then, $\partial (\phi + \psi) = \overline{{\mathcal A} + {\mathcal B}}$. As a consequence (1.1) has a unique solution.
\end{remark}
Let ${\mathcal A}$ be selfadjoint monotone operator and
let ${\mathcal B}$ be a nonlinear maximal monotone operator (possibly multi-valued)
on $(L^2(0 , T; {\mathbb H}) \; ; \; \langle \langle \; , \; \rangle \rangle)$. They said forming an acute angle if the following
holds:
$$(2.7) \; \; \; \; \;\; (( {\mathcal A_{\lambda}} u , {\mathcal B_{\mu}}
u \; )) \; \geq \; 0, \; \; \forall \lambda, \mu > 0, \; \; u \in 
 L^2(0 , T; {\mathbb H})$$
 
\begin{theorem}
{\it Let ${\mathcal A}$ be a linear self-adjoint monotone operator and
let ${\mathcal B}$ be a nonlinear maximal monotone operator (possibly multi-valued). Assume that
(2.7) holds. Then $({\mathcal A} + {\mathcal B})_{v} \equiv {\mathcal A} + {\mathcal B}$ is an m-accretive nonlinear operator on $L^2(0 , T; {\mathbb H})$.}
\end{theorem}
\begin{proof}
Since ${\mathcal A} + {\mathcal B_{\mu}}$ is m-accretive, then for all
$w \in L^2(0 , T; {\mathbb H})$, there exists a unique $(u_{\mu})_{\mu >
0} \in D({\mathcal A})$ such that

$$(2.8) \; \; \; \; \; \; \mu u + {\mathcal A} u_{\mu} + {\mathcal B_{\mu}}u_{\mu}
= w, \; \; \mu > 0$$
According to Br\'ezis-Grandall-Pazy (see \cite{br-gr-pa}), the problem
$$(2.9) \; \; \; \; \; \; u + {\mathcal A} u + {\mathcal B}u
\ni w$$
has a unique solution $u \in D({\mathcal A}) \cap D({\mathcal B})$ if and only
if the family $({\mathcal B_{\mu}})_{\mu > 0}$ is bounded in 
$L^2(0 , T; {\mathbb H})$. Now, a sufficient condition for $({\mathcal B_{\mu}})_{\mu > 0}$ to be bounded is guaranteed by (2.7).

Since $({\mathcal B_{\mu}})_{\mu > 0}$ is bounded, this implies that $u_{\mu}$ strongly
converges to some $ u \in L^2(0 , T; {\mathbb H})$ and that
${\mathcal B_{\mu}} u_{\mu}$ weakly converges to $Z$ as $\mu$ goes to $0$.
Using the fact ${\mathcal B}$ is a maximal monotone operator, it easily follows
that $u \in D({\mathcal B})$ and that ${\mathcal B} u = Z$. Using a
similar argument, it turns out that $u \in D({\mathcal A})$ and that
$u$ satisfies (2.9), that is $({\mathcal A} + {\mathcal B})_{v} \equiv {\mathcal A} + {\mathcal B}$ is a nonlinear maximal monotone operator on $L^2(0 , T; {\mathbb H})$
\end{proof}

\begin{corollary}
{\it Let ${\mathcal A}$ be a linear self-adjoint monotone operator and 
let ${\mathcal B}$ be a nonlinear maximal monotone operator (possibly multi-valued). Assume that
(2.7) holds. Then the problem (1.1) has a unique solution.}
\end{corollary}

\subsection{Examples}
\subsubsection{Example 1}
Let
$\Omega \subset {\mathbb R^n}$ be a bounded open subset with smooth boundary
$\partial \Omega$.  Consider the partial differential equations of type
\[ (2.10) \; \; \; \; \; \; \; \left \{ \begin{array}{cl}
u_{t} + A u (t,x) + B u(t,x)  \ni f(t,x), \; \; \forall (t,x) \in (0, T] \times \Omega\\
u (0 , x) = 0, \; \; \forall x \in \partial \Omega
\end{array} \right.
\]
where $A = - \Delta$ and $B u = F(u)$ with the normalization $F(0) \ni 0$
in $L^2(\Omega)$. The function $f$ belongs to $L^2(0 , T; 
L^2(\Omega))$, where $L^2(0 , T;
L^2(\Omega))$ is the Hilbert space endowed with the inner product
$$ \langle \langle u \; , \; v \rangle \rangle = \int_{0}^{T} \; \langle
u(t,x) \; , \; v(t,x) \rangle_{L^2(\Omega)}
dt$$
It is convenient to write (2.10) of the form
$$(2.11) \; \; \; \; \; \; \; {\mathcal S} u + {\mathcal A} u + {\mathcal B} u \ni f$$
where
\[ \left \{ \begin{array}{cl}
D({\mathcal S}) = \{ u \in {\mathbb H^1}(0 , T; L^2(\Omega)) \; | \; u(0, x) = 0, \; \forall x \in \partial \Omega \}\\
{\mathcal S} u = u_{t}, \; \; \; \forall u \in D({\mathcal S})
\end{array} \right.
\]
Recall that $A$ and $B$ are 
respectively given in the Hilbert space $L^2(\Omega)$ by
$$D(A) = {\mathbb H_{0}^{1}}(\Omega) \cap {\mathbb H^{2}}(\Omega)
\; \; \; \mbox{and} \; \; \; A u = - \Delta u, \; \; \forall
u \in D(A)$$
and
\[\left \{ \begin{array}{cl}
D(B) = \{  u \in L^2(\Omega)) \; | \; F(u)
\in L^2(\Omega))\}\\
B u(x) = F(u)(x) \; \; \mbox{a.e} \; \; u \in D(B)
\end{array} \right.
\]
and that ${\mathcal A}$ and ${\mathcal B}$ are 
defined in the Hilbert spce $L^2(0 , T; L^2(\Omega))$ as described in the introduction.

We assume that $F: \; {\mathbb R} \mapsto {\mathbb R}$ is everywhere defined
non-decreasing function of class $C^1$ satisfying the assumption $F(0) =
0$. Under such a assumption, then the operator $B$ is m-accretive on
$L^2(\Omega)$ (see  [ Proposition 2.7, \cite{br}]. Thus, the corresponding
operators ${\mathcal A}$ and ${\mathcal B}$ are respectively m-accretive
operators on $L^2( 0 , T; L^2(\Omega))$.

We have

\begin{proposition}
{\it 
Let ${\mathcal A}$ and ${\mathcal B}$ be the operators described above.
Assume that $F: \; {\mathbb R} \mapsto {\mathbb R}$ is everywhere defined
non-decreasing function of class $C^1$ satisfying the assumption $F(0) =
0$ and that $D({\mathcal A}) \cap D({\mathcal B}) \not= \emptyset$. Then
the variational sum of ${\mathcal A}$ and ${\mathcal B}$ is a maximal monotone
operator and that
$$( {\mathcal A} + {\mathcal B})_{v} =  {\mathcal A} + {\mathcal B}$$
}

\end{proposition}

\begin{proof}
It is sufficient to prove that the problem given as
\[(2.12) \; \; \; \; \; \; \; \left \{ \begin{array}{cl}
\mu u - \Delta u + F(u) = f  \; \; \in \Omega, \; \; \mu > 0\\
u = 0 \; \; \; \; \; \mbox{on} \; \; \partial \Omega
\end{array} \right.
\]
has a unique solution for every $f \in L^2(\Omega)$. 

The existence and the uniqueness of a solution to (2.12) is guaranteed by
a result of Br\'ezis - Grandall - Pazy (see [Theorem 3.1, \cite{br-gr-pa} ]), which says that,
under previous assumptions one has 
$$\int_{\Omega} - \Delta u (x) B_{\lambda} u (x) dx \; \geq \; 0,
\; \; \forall u \in  {\mathbb H_{0}^{1}}(\Omega) \cap {\mathbb H^{2}}(\Omega),
\; \; \lambda > 0$$
That is, $- \Delta + B$ is a maximal monotone operator. Therefore
${\mathcal A} + {\mathcal B}$ is also a maximal monotone operator. According to the definition of the
variational sum of ${\mathcal A}$ and ${\mathcal B}$, clearly
$({\mathcal A} + {\mathcal B})_{v} = {\mathcal A} + {\mathcal B}$ is a maximal
monotone operator.
\end{proof}

\begin{corollary}
{\it 
Let ${\mathcal A}$ and ${\mathcal B}$ be the operators described above.
Assume that $F: \; {\mathbb R} \mapsto {\mathbb R}$ is everywhere defined
non-decreasing function of class $C^1$ satisfying the assumption $F(0) =
0$ and that $D({\mathcal A}) \cap D({\mathcal B}) \not= \emptyset$. Then
the problem (2.10) has a unique solution.}
\end{corollary}
\begin{proof}
The problem (2.10) is equivalent to the problem (2.11). According to
Proposition 2.8, the solvability of (2.11) is established.
\end{proof}
\subsubsection{Example 2}
Consider the problem (2.10) where both $A$ and $B$ are linear operators on
${\mathbb H} = L^2({\mathbb R^n})$ defined as
$$D(A) = {\mathbb H^{2}}({\mathbb R^n})
\; \; \; \mbox{and} \; \; \; A u = - \Delta u, \; \; \forall
u \in D(A)$$
and
\[\left \{ \begin{array}{cl}
D(B) = \{  u \in L^2({\mathbb R^n}) \; \; | \; \; Q(x) u
\in L^2({\mathbb R^n})\}\\
B u = Q(x) u \; \; u \in D(B)
\end{array} \right.
\]
where the potential $Q$ satisfies the following assumptions
$$(2.13) \; \; \; \; \; \; Q(x) > 0, \; \; \; Q \in L^1({\mathbb
R^n}), \; \; \mbox{and} \; \; Q \not \in L_{loc}^2({\mathbb
R^n})$$

\begin{proposition}
Under assumption (2.13), then $D(A) \cap D(B) = \{ 0 \}$.
\end{proposition}

\begin{proof}
The proof of this proposition depends on the dimensional space
$n$ (This is explained by the Sobolev embedding ) . We will provide
the proof in the case where $n \leq 3$. Indeed the proof in the case where
$n \geq 4$ can be found in [ Proposition 2.1, \cite{dia4} ].

Let $u \in D( - \Delta ) \cap D( Q )$ and assume that $u \not \equiv 0$.
Since $u \in {\mathbb H^2}({\mathbb R^n})$ where $n \leq 3$, then
$u$ is a continuous function by Sobolev theorem (see \cite{adam}).
There exists an opren subset $\Omega$ of ${\mathbb R^n}$ and there exists
$\delta > 0$ such that  $| u(x) | >  \delta$ for all $x \in \Omega$.
Let $\Omega'$ be a compact subset of $\Omega$, equipped with the induced
topology by $\Omega$ ($\Omega'$ is a compact subset of ${\mathbb R^n}$).
It easily follows that
$$ | Q |_{|_{\Omega'}} = {\frac{ (| Q u |)_{|_{\Omega'}}}{| u |_{|_{\Omega'}}}}
\in L^2(\Omega')$$
since $(| Q u |)_{|_{\Omega'}} \in L^2(\Omega')$ and 
${\frac{ 1 }{| u |_{|_{\Omega'}}}} \in L^{\infty}(\Omega')$. Therefore $Q
\in L^2(\Omega')$; this is impossible according to the assumption (2.13),
then $u \equiv 0$.
\end{proof}

\noindent {\bf Example of potential $Q$ satisfying (2.13)}.
Let $\Omega$ be a compact subset of ${\mathbb R^n}$ and let $G$ be a complex
function satisfying, $ \Re e \; G > 0, \; G \in L^1(\Omega), \; G \not \in
L^2(\Omega) \; \mbox{and} \; G \equiv 0 \; \mbox{on} \; {\mathbb R^n} \; - \; \Omega$.
Consider the following rational sequence $\alpha_{k} = ( \alpha_{k}^{1}, \alpha_{K
}^{2}, ......, \alpha_{k}^{n}) \in {\mathbb Q^n}$. Then the function $Q$ given by,

$$Q(x) = \sum_{k = 1}^{+ \infty} \frac{G( x - \alpha_{k} )}{k^2},$$
satisfies the assumption (2.13).

\bigskip 
 
Note that $A$ and $B$ given above are respectively self-adjoint operators
on $L^2({\mathbb R^n})$. In what follows, we shall assume that
$D(A^{\frac{1}{2}}) \cap D(B^{\frac{1}{2}})$ is dense in $L^2({\mathbb R^n})$;
one can show that the variational sum $(A+B)_{v}$ and the sum form $A \oplus
B$ of
$A$ and $B$ coincide, that is,
$$(A + B )_{v} \equiv  A \oplus B$$

Consider the closed sesquilinear forms given by
$$\Phi(u,v) = \int_{R^n} \nabla u \; \bar{\nabla v} dx \; \; \mbox{for all}
\; \; u,v \; \in {\mathbb H^1}({\mathbb R^n}),$$
$$\Psi(u,v) = \int_{\mathbb R^n} Qu \bar{v} dx \; \; \mbox{for all} \; u,v \; \in D(B^{\frac{1}{2}}),$$
and the sum of the forms $\Phi$ and $\Psi$ is given by, $\Upsilon = \Phi + \Psi$, in other words,

$\Upsilon(u,v) = \int_{\mathbb R^n} ( \nabla u \bar{\nabla v} + Qu \bar{v})dx \; \;
\mbox{for all} \; u,v \; \in D(A^{\frac{1}{2}}) \cap D(B^{\frac{1}{2}}).$
The sesquilinear form $\Upsilon$ is a closed sectorial and densely defined form as sum of 
closed sectorial and densely defined forms, then there exists a unique m-sectorial
operator $A \oplus B$, called sum form of $A$ and $B$
associated to $\Upsilon$ ( see \cite{dia1} ) and $\Upsilon$ has the following represention,
$$\Upsilon(u,v) = \langle (A \oplus B)u \; , \; v \rangle\; \; \mbox{for all}
\; u \; \in D(A \oplus B), \; v \; \in D(A^{\frac{1}{2}}) \cap D(B^{\frac{1}{2}})$$
According to the author ( \cite{dia1} ), the operator $A \oplus B$ verifies the
well-known condition of Kato, in other words,
$$ D((A \oplus B)^{\frac{1}{2}}) = D(\Upsilon) = D(((A \oplus B)^*)^{\frac{1}{2}})$$
The operator $A \oplus B$ has been computed by H. Br\'ezis and T. Kato in  \cite{brezis-kato}. It
is given by
\[\left \{ \begin{array}{cl}
D(A \oplus B) = \{ u \in {\mathbb H^1}({\mathbb R^n}) \; \; | \; \; Q|u|^2 \in L^1({\mathbb R^n}), \; \; \;
- \Delta u + Qu \in L^2({\mathbb R^n}) \} \\
 (A \oplus B)u = - \Delta u + Qu
\end{array} \right.
\]
Clearly $(A + B)_{v} \equiv ( A \oplus B )$. 

\bigskip

Now, since $D(A) \cap D(B) = \{ 0 \}$ under (2.13), the problem (2.10) does
not make sense anymore. Another alternative is to consider the following
problem
\[ (2.14) \; \; \; \; \; \; \; \left \{ \begin{array}{cl}
u_{t} + (A \oplus B)u(t,x)  = f(t,x), \; \; \forall (t,x) \in (0, T] \times \Omega\\
u (0 , x) = 0, \; \; \forall x \in \partial \Omega
\end{array} \right.
\]
Clearly, the problem (2.14) has a unique solution.

\bigskip
Let us define ${\mathcal A} \oplus {\mathcal B} = ({\mathcal A} + {\mathcal
B})_{v}$. Let
${\mathcal A}$ and ${\mathcal B}$ be the corresponding operators described
in the introduction and defined in the Hilbert space $L^2( 0 , T; L^2({\mathbb
R^n}))$. The corresponding operator to the sum form ${\mathcal A} \oplus {\mathcal B}$ is defined
on $L^2( 0 , T; L^2({\mathbb
R^n}))$ by:
$u \in D({\mathcal A} \oplus {\mathcal B})$ iff $u \in L^2( 0 , T; {\mathbb H^1}({\mathbb R^n})), \; \; Q|u|^2 \in L^2( 0 , T; L^1({\mathbb R^n})), \; \; \mbox{and} \; \;
- \Delta u + Qu \in L^2( 0 , T; L^2({\mathbb R^n}))$ with
$({\mathcal A} \oplus {\mathcal B})u = - \Delta u + Qu$ $\blacksquare$

\bibliographystyle{amsplain}

\end{document}